\documentclass[12pt]{amsart}

\usepackage{amsthm,xypic,amssymb,verbatim}

\numberwithin{equation}{section}
 \newtheorem{thm}{Theorem}[section]
 \newtheorem{cor}[thm]{Corollary}
 \newtheorem{lem}[thm]{Lemma}
 
\newtheorem{problem}[thm]{Question}
 \theoremstyle{definition}
 
 \theoremstyle{remark}
 \newtheorem{rem}[thm]{Remark}
 \theoremstyle{definition}
 \newtheorem{ex}[thm]{Example}

 \newtheorem*{thmNONUMBER}{Theorem}

 \newcommand{\CC}{\mathbb{C}}

 \newcommand{\XX}{\mathbb{X}}
 \newcommand{\PP}{\mathbb{P}}

\begin{document}

\title{Star configuration points and generic plane curves}

\author[E. Carlini]{Enrico Carlini}
\address[E. Carlini]{Dipartimento di Matematica, Politecnico di Torino, Turin, Italy}
\email{enrico.carlini@polito.it}

\author[A. Van Tuyl]{Adam Van Tuyl}
\address[A. Van Tuyl]{Department of Mathematical Sciences,
Lakehead University, Thunder Bay, ON, Canada, P7B 5E1}
\email{avantuyl@lakeheadu.ca}

%%% ----------------------------------------------------------------------

\begin{abstract}
Let $\ell_1,\ldots\ell_l$ be $l$ lines in $\PP^2$ such that no
three lines meet in a point.  Let $\XX(l)$ be the set of points
$\{\ell_i \cap \ell_j ~|~ 1 \leq i < j \leq l\} \subseteq \PP^2$.
We call $\XX(l)$ a star configuration.  We describe
all
pairs
$(d,l)$ such that the generic
degree $d$ curve in $\PP^2$ contains a
$\XX(l)$.
\end{abstract}

%%% ----------------------------------------------------------------------
\maketitle
%%% --------------------------------------------------------------------

\section{Introduction}

The problem of studying subvarieties of algebraic varieties is a
crucial one in algebraic geometry, e.g., the case
of divisors.  The study of subvarieties of hypersurfaces in $\PP^n$ has a
particularly rich  history. For example, one can look for the existence
of $m$ dimensional linear spaces on generic hypersurfaces of
degree $d$ leading to the theory of Fano varieties and to the well
known formula relating $n,m$ and $d$ (e.g., see \cite[Theorem
12.8]{Harris}).

Because linear spaces are complete intersections, it is natural to
look for the existence of complete intersection subvarieties on a
generic hypersurface. The case of codimension two complete
intersections was first studied by by Severi \cite{Sev06}, and
later generalized and extended to higher codimensions by Noether,
Lefschetz \cite{Lef21} and Groethendieck \cite{Gro05}.  Recently,
in \cite{CCG1}, secant vavieties
were used to give a complete solution for the existence of complete
intersections of codimension $r$ on generic hypersurfaces in
$\PP^n$ when $2r\leq n+2$.
Fewer results are known if the codimension of the
complete intersection is large, i.e., the codimension is close to
the dimension of the ambient space. In \cite{Szabo} the case of
complete intersection curves is studied and completely solved. The
case of complete intersection points on generic surfaces in
$\PP^3$ is considered in \cite{CCG2} where some asymptotic results
are presented. The case of complete intersection points in $\PP^2$
is a special case of \cite{CCG1}, and it is completely solved.

Taking our inspiration from \cite{CCG1,CCG2}, we examine the
problem of determining when special configurations of points lie
on a generic degree $d$ plane curve in $\PP^2$.  We shall focus on
{\bf star configurations} of points. Consider $l$ lines in the
plane $\ell_1,\ldots,\ell_l\subset\PP^2$ such that
$\ell_i\cap\ell_j\cap\ell_k=\emptyset$. The set of points
$\mathbb{X}(l)$ consisting of the ${l\choose 2}$ pairwise
intersections of the lines $\ell_i$ is called a star
configuration. In this paper, we address the following question

\begin{problem}\label{mainproblem}
For what pairs $(d,l)$ does the generic degree $d$ plane curve
contain a star configuration $\mathbb{X}(l) \subseteq \PP^2$?
\end{problem}

The name star configuration was suggested by A.V. Geramita
since $\XX(5)$ is the ten points of intersections when
drawing a star using five lines (see Figure 6.2.4 in
\cite{CHT}).  The configurations
$\XX(l)$ appeared in the work of Geramita, Migliore, and Sabourin
\cite{GMS} as the support of a set of double points
whose Hilbert function exhibited an extremal behaviour.
More recently, Cooper, Harbourne, and Teitler \cite{CHT} computed
the Hilbert function of any homogeneous set of fatpoints supported on
$\XX(l)$.  Bocci and Harbourne \cite{BH} used star
configurations (and their generalizations) to compare
the symbolic and regular powers of an ideal.  Further properties
of star configurations continue to be uncovered;  e.g., ongoing work
of Geramita, Harbourne,
Migliore \cite{GHM}.

We can answer Question \ref{mainproblem} because we can exploit the
rich algebraic structure of star configurations.  In particular,
we will require the fact that one can easily write down a list of
generators for $I_{\XX(l)}$, the defining ideal of $\XX(l)$,
as well as the fact that the Hilbert function of
$\XX(l)$ is the same as the Hilbert function of $\binom{l}{2}$
generic points in $\PP^2$. Using these properties, among others,
we give the following solution to Question \ref{mainproblem}:

\begin{thmNONUMBER}{\bf{\ref{maintheorem}} }{\it
Let $l \geq 2$.  Then the generic degree $d$ plane curve
contains a star configuration $\XX(l)$ if and only if
\begin{enumerate}
\item[$(i)$] $l=2$ and $d \geq 1$, or
\item[$(ii)$] $l=3$ and $d \geq 2$, or
\item[$(iii)$] $l=4$ and $d \geq 3$, or
\item[$(iv)$] $l = 5$ and $d \geq 5$.
\end{enumerate}}
\end{thmNONUMBER}

Our proof is broken down into a number of cases.  It will
follow directly from the generators of
$I_{\XX(l)}$ that Question \ref{mainproblem} can have no
solution for $d < l-1$.  Using a simple dimension
counting argument, Theorem \ref{asymptoticresult} shows that for $l \geq 6$,
there is no solution to Question \ref{mainproblem}.  The
cases $l=2$ and $l=3$ are trivial cases, so the bulk of the
paper will be devoted to the cases that $l=4$ and $l=5$.

To prove these cases, we rephrase Question \ref{mainproblem}
into a purely ideal theoretic question (see Lemma \ref{tgspaceideallem}).
Precisely, we construct a new ideal $I$ from the generators of $I_{\XX(l)}$
and the linear forms defining the lines $\ell_1,\ldots,\ell_l$.
We then show that Question \ref{mainproblem} is equivalent
to determining whether $I_d = (\CC[x,y,z])_d$.  We then
answer this new algebraic reformulation.

The proofs for the cases $(d,l) = (3,4)$ and $(4,5)$ are, we
believe, especially interesting. To prove that Question
\ref{mainproblem} is true for $(3,4)$, we exploit the natural
group structure of the plane cubic curve to find a $\XX(4)$ on the
curve.  To show the non-existence of a solution for
$(4,5)$, we require the classical theory of {\bf L\"uroth
quartics}. L\"uroth quartics are the plane quartics that pass
through a $\XX(5)$. L\"uroth quartics have the property of forming a
hypersurface in the space of plane quartics.  The existence of
this hypersurface is the obstruction for the existence of a
solution when $(d,l) = (4,5)$.

Our paper is structured as follows.  In Section 2,
we describe the needed algebraic properties of star configurations.  In
Section 3, we give some asymptotic results.  In Section 4, in
preparation for the last two sections, we rephrase Question
\ref{mainproblem} into an equivalent algebraic problem.  Sections 5
and 6 deal with the cases $l=4$ and $l=5$, respectively.

\noindent {\bf Acknowledgements}  This paper began when the first
author visited the second at Lakehead University.  Theorem
\ref{maintheorem} was inspired by computer experiments using CoCoA
\cite{C}. The first author was partially supported by the Giovani
Ricercatori grant 2008 of the Politecnico di Torino. Both authors
acknowledge the financial support of NSERC.

\section{Star Configurations}

Throughout this paper, we set $S = \CC[x,y,z]$ and we denote its
$d$th homogeneous piece with $S_d$. Moreover, we fix standard
monomial bases in each $S_d$ in such a way that $\PP S_d$ will be
identified with $\PP^{N_d}$ where $N_d={d+2 \choose 2}-1$. We
recall the relevant definitions concerning star configurations of
points in $\PP^2$.

Let $l\geq 2$ be an integer. The scheme
$\mathbb{X}(l)\subset\PP^2$ is said to be a {\bf star
configuration} if it consists of ${l\choose 2}$ distinct points
which are the pairwise intersections of $l$ distinct lines, say
$\ell_1,\ldots,\ell_l$, where no three lines pass through the same
point.   We will also call $\mathbb{X}(l)$ a {\bf star
configuration set of points.}

Note that when $l=2$, then $\XX(2)$ consists of a single point.
When $l=3$, then $\XX(3)$ is any set of three points in $\PP^2$,
provided the three points do not lie on the same line.  It is
clear that any degree $d \geq 1$ plane curve contains a point.
Furthermore, when $d \geq 2$, the generic degree $d$ plane curve
will contain three points not lying on a line.  These remarks take
care of the trivial cases of Theorem \ref{maintheorem}, which we
summarize as a lemma:

\begin{lem}  \label{trivialcase}
The generic degree $d$ plane curve contains a star configuration $\XX(l)$
with $l =2$, respectively $l=3$, if and only if $d \geq 1$, respectively $d \geq 2$.
\end{lem}

Given a $\XX(l)$, for each $i=1,\ldots,l$, we let $L_i$ denote a linear form in
$S_1$ defining the line $\ell_i$. The defining
ideal of $I_{\XX(l)}$ is then given by
\[I_{\mathbb{X}(l)}=\bigcap_{i\neq j} (L_i,L_j).\]
We will sometimes write a point of $\XX(l)$ as $p_{i,j}$, where $p_{i,j}$
is the point defined by the ideal $(L_i,L_j)$, i.e., $p_{i,j}$ is
the point of intersection of the lines $\ell_i$ and $\ell_j$.
The following lemma describes the generators
of $I_{\XX(l)}$;  we will exploit this fact throughout the paper.

\begin{lem}  Let $l \geq 2$, and let $\XX(l)$ denote
the star configuration constructed from the lines
$\ell_1,\ldots,\ell_l$.  If $L_i$ is a linear form defining
$\ell_i$ for $i = 1,\ldots,l$, then
\[I_{\XX(l)} = (\hat L_1,\ldots,\hat L_l) ~~\mbox{where}~~ \hat L_i=\prod_{j\neq i} L_j.\]
\end{lem}

For a proof of this fact, see \cite[Claim in Proposition
3.4]{miochow}.  From this description of the generators, we see
that $I_{\XX(l)}$ is generated in degree $l-1$. Thus, for $d <
l-1$, there are no plane curves of degree $d$ that contain a
$\XX(l)$ since $(I_{\XX(l)})_d = (0)$. Hence, as a corollary, we
get some partial information about Question \ref{mainproblem}.

\begin{cor}\label{d<l-1rem}
If $d < l-1$, then the generic degree $d$ curve does not contain a $\XX(l)$.
\end{cor}

It is useful to recall that star configuration points have the
same Hilbert function as generic points.
More precisely, by \cite[Lemma 7.8]{GMS}, we have

\begin{lem}\label{hf}
Let $HF(\mathbb{X}(l),t) = \dim_\CC (S/I_{\XX(l)})_t$ denote the
Hilbert function of $S/I_{\XX(l)}$.  Then
\[HF(\mathbb{X}(l),t)=\min\left\{{t+2\choose 2}, {l\choose 2}\right\}
~~\mbox{for all $t \geq 0$}.\]
\end{lem}

\section{An asymptotic result}

In \cite{CCG1} it is shown that the generic
degree $d$ plane curve contains a $0$-dimensional complete intersection
scheme of type $(a,b)$ whenever $a,b\leq d$
(actually the result holds for {\em any} degree $d$ plane curve).
Thus, arbitrarily large complete intersections can be found on
generic plane curves of degree high enough. But the same does not hold for star configurations
as shown by Theorem \ref{asymptoticresult}.

Using the previous description of $I_{\mathbb{X}(l)}$ we introduce
a quasi-projective variety parametrizing star configurations.
Namely, we consider
\[\mathcal{D}_l\subset \underbrace{\check\PP^2\times\cdots\times \check\PP^2}_{l}\]
such that $(\ell_1,\ldots,\ell_l)\in \mathcal{D}_l$ if and only if
no three of the lines $\ell_i$ are  passing through the same
point; here, $\check\PP^2$ denotes the dual projective space.

\begin{thm}\label{asymptoticresult}
If $l>5$ is an integer, then the generic degree $d$ plane curve
does not contain any star configuration $\mathbb{X}(l)$.
\end{thm}

\begin{proof}
It is enough to consider the case $d\geq l-1$ (see Corollary
\ref{d<l-1rem}). Let $\PP S_d$ be the space parametrizing degree
$d$ planes curves and define the following incidence
correspondence
\[\Sigma_{d,l}=\left\{(\mathcal{C},\mathbb{X}(l)) :\mathcal{C}\supset\mathbb{X}(l)\right\}\subset\PP S_d\times\mathcal{D}_l.\]
We also consider the natural projection maps
\[\psi_{d,l}:\Sigma_{d,l}\longrightarrow \mathcal{D}_l
~~\text{and}~~~
\phi_{d,l}:\Sigma_{d,l}\longrightarrow \PP S_d.\]
Clearly we have that $\phi_{d,l}$ is dominant if and only if
Question \ref{mainproblem} has an affirmative answer.

Using a standard fiber dimension argument, we see that
\[\dim\Sigma_{d,l}\leq\dim\mathcal{D}_l+\dim\psi_{d,l}^{-1}(\mathbb{X}(l)) \]
for a generic star configuration $\mathbb{X}(l)$, where
$\dim\mathcal{D}_l = 2l$
and by Lemma \ref{hf}
\[\dim\psi_{d,l}^{-1}(\mathbb{X}(l))=\dim_\CC
\left(I_{\XX(l)}\right)_d-1={d+2 \choose 2}-{l \choose 2}-1.\]

Hence the map $\phi_{d,l}$ is dominant only if $\dim\Sigma_{d,l}-\dim\PP
S_d\geq 0$ and this is equivalent to
\[2l-{l \choose 2}=\frac{l(5-l)}{2}\geq 0.\]
Thus the result is proved.
\end{proof}

By Lemma \ref{trivialcase} and the above result,
we only need to treat the cases $l=4$ and $5$.  We postpone
these cases to first rephrase Question \ref{mainproblem}
into an equivalent algebraic question.

\section{Restatement of Question \ref{mainproblem}}

We derive some technical results, moving our
Question \ref{mainproblem}  back and forth between questions in algebra and
questions in geometry.  First, we notice the following trivial fact:

\begin{lem}\label{fsummililem} Let $\{F=0\}$ be an equation of the degree $d$ curve
$\mathcal{C}\subset\PP^2$. Then $\mathcal{C}$ contains a star
configuration $\mathbb{X}(l)$ only if
\[F=\sum_{i=1}^l M_i\hat L_i\]
where the forms $M_i$ have degree $d-l+1$ and the forms $\hat L_i$
are defined as $\hat L_i=\prod_{j\neq i}L_i$ for some linear forms
$L_1,\ldots,L_l$.
\end{lem}

Hence, it is natural to perform the following geometric
construction. We define a map of affine varieties
\[\Phi_{d,l}:\underbrace{ S_1\times\cdots\times S_1}_{l}\times\underbrace{ S_{d-l+1}\times\cdots\times S_{d-l+1}}_{l}\longrightarrow S_d\]
such that
\[\Phi_{d,l}\left(L_1,\ldots,L_l,M_1,\ldots,M_l\right)=\sum_{i=1}^l M_i\hat
L_i\] We then rephrase our problem in terms of the map
$\Phi_{d,l}$:

\begin{lem}\label{dominantmaplemma} Let $d,l$ be nonnegative integers
with $d \geq l-1$. Then the following
are equivalent:
\begin{itemize}
\item[$(i)$] Question \ref{mainproblem} has an affirmative answer for $d$ and
$l$;
\item[$(ii)$] the map $\Phi_{d,l}$ is a dominant map.
\end{itemize}
\end{lem}
\begin{proof}
Lemma \ref{fsummililem} proves that $(i)$ implies $(ii)$. To prove
the other direction, it is enough to show that for a generic form
$F$, the fiber $\Phi_{d,l}^{-1}(F)$ contains  a set of $l$ linear
forms defining a star configuration. More precisely, define
$\Delta\subset  S_1\times\cdots\times S_1\times
S_{d-l+1}\times\cdots\times S_{d-l+1}$ as follows:
\[\Delta = \left\{\left(L_1,\ldots,L_l,M_1,\ldots,M_l\right)~\left|
\begin{tabular}{l}
\mbox{there exists $a \neq b \neq c$ such that} \\
\mbox{$L_a,L_b,L_c$ are linearly dependent}
\end{tabular} \right\}\right. .\]
Then we want to show that $\Phi_{d,l}^{-1}(F)\not\subset\Delta$.

We proceed by contradiction, assuming that the generic fiber of
$\Phi_{d,l}$ is contained in $\Delta$.  Then $\Delta$ would be a
component of the domain of $\Phi_{d,l}$. Thus a contradiction as
the latter is an irreducible variety being the product of
irreducible varieties.
\end{proof}

Using the map $\Phi_{d,l}$ we can now translate Question \ref{mainproblem} into an
ideal theoretic question.

\begin{lem}\label{tgspaceideallem}
Let $d,l$ be non-negative integers. Consider generic forms
$L_1,\ldots,L_l\in S_1$ and $M_1,\ldots,M_l\in S_{d-l+1}$. Set
\[\hat L_i=\prod_{j\neq i}L_j
~~\mbox{and}~~
\hat L_{i,j}=\prod_{h\neq i,h\neq j}L_h, \mbox{ for } i\neq j.\]
Also set
\begin{eqnarray*}
Q_1 &= &M_2\hat L_{1,2}+\cdots+M_l\hat L_{1,l}=\sum_{i\neq 1} M_i \hat L_{1,i},\\
&\vdots& \\
Q_l&=&M_1\hat L_{l,1}+\cdots+M_{l-1}\hat L_{l,l-1}=\sum_{i\neq l} M_i \hat L_{l,i}.
\end{eqnarray*}
With this notation, form the ideal
\[I=(\hat L_1,\ldots,\hat L_l,Q_1,\ldots,Q_l)\subset S.\]
Then the following are equivalent:
\begin{itemize}
\item[$(i)$] Question \ref{mainproblem} has an affirmative answer for $d$ and
$l$;
\item[$(ii)$] $I_d=S_d$.

\end{itemize}
\end{lem}
\begin{proof}
Using Lemma \ref{dominantmaplemma} we just have to show that
$\Phi_{d,l}$ is a dominant map. In order to do this we will
determine the tangent space to the image of $\Phi_{d,l}$ in a
generic point $q=\Phi_{d,l}(p)$, where
$p=(L_1,\ldots,L_l,M_1,\ldots,M_l)$ and we denote with $T_q$ this
affine tangent space.

The elements of the tangent space $T_q$ are obtained as

\[\left.{d \over dt}\right|_{t=0} \Phi_{d,l}\left(L_1+tL_1',\ldots,L_l+tL_l',M_1+tM_1',\ldots,M_l+tM_l'\right)\]
when we vary the forms $L_i'\in S_1$ and $M_i'\in S_{d-l+1}$. By a
direct computation we see that the elements of $T_q$ have the form
\[M_1'\hat L_1+\cdots+M_l'\hat L_l+ L_1'(M_2\hat L_{1,2}+\cdots+M_l\hat L_{1,l}) +\cdots+ \]
\[+L_j'(M_1\hat L_{1,j}+\cdots+M_l\hat L_{j,l})+\cdots+L_l'(M_1\hat L_{1,2}+\cdots+M_{l-1}\hat L_{1,l-1}),\]
where $\hat L_i=\prod_{j\neq i}L_i$ and $\hat L_{i,j}=\prod_{h\neq
i,h\neq j}L_h, \mbox{ for } i\neq j$.

Since the $L_i'\in S_1$ and $M_i'\in S_{d-l+1}$ can be
chosen freely, we obtain that $I_d=T_q$.
\end{proof}

\begin{rem}\label{cocoacomputREMARK}
Lemma \ref{tgspaceideallem} is an effective tool to give a
positive answer to each special issue of Question
\ref{mainproblem}. Given $d$ and $l$ we construct the ideal $I$ as
described by choosing forms $L_i$ and $M_i$. Then we compute $\dim_\CC
I_d$ using a computer algebra system, e.g., CoCoA. If $\dim_\CC
I_d=\dim_\CC S_d$, by upper semi-continuity of the dimension, we have
proved that our question has an affirmative answer for these given
$d$ and $l$. In fact, the dimension can decrease only on a proper
closed subset. On the contrary, if $\dim_\CC I_d<\dim_\CC S_d$ we do not
have any proof. Question \ref{mainproblem} can both have a negative (for any
choice of forms the inequality holds) or an affirmative answer
(our choice of forms is too special and another choice will give
an equality). But, if we choose our forms generic enough, we do
have a strong indication that the answer to Question \ref{mainproblem} should be
negative.
\end{rem}

\begin{ex}\label{l4d4ex}
We use CoCoA \cite{C} to
consider the example $(d,l) = (4,4)$.
\footnotesize
\begin{verbatim}
--we define the ring we want to use
Use S::=QQ[x,y,z];
--we choose our forms
L1:=x;
L2:=y;
L3:=z;
L4:=x+y+z;
M1:=x+y-z;
M2:=-x+2y+2z;
M3:=2x-y-z;
M4:=x+y+2z;
--we build the forms Qi
Q1:=M2*L3*L4+M3*L2*L4+M4*L2*L3;
Q2:=M1*L3*L4+M3*L1*L4+M4*L1*L3;
Q3:=M1*L3*L4+M2*L1*L4+M4*L1*L2;
Q4:=M1*L2*L3+M2*L1*L3+M3*L1*L2;
--we define the ideal I
I:=Ideal(L2*L3*L4,L1*L3*L4,L1*L2*L4,L1*L2*L3,Q1,Q2,Q3,Q4);
--we compute the Hilbert function of S/I in degree 4
Hilbert(S/I,4);
\end{verbatim}
\normalsize
Evaluating the code above we get
\footnotesize
\begin{verbatim}
Hilbert(S/I,4);
0
-------------------------------
\end{verbatim}
\normalsize
and this means that $\dim_\CC S_4-\dim_\CC I_4=0$. Hence we showed that the generic plane quartic contains a $\XX(4)$.
\end{ex}

\section{The $l=4$ case}

In this  section we consider Question \ref{mainproblem} when $l
=4$. We begin with a special instance of Question \ref{mainproblem}, that is,
when $d=3$, since it has a nice geometric proof which takes
advantage of the group structure on the curve.

\begin{lem}\label{grouplawlem}
Let $\mathcal{C}\subset\PP^2$ be a smooth cubic curve. Then there
exists a star configuration $\mathbb{X}(4)\subset \mathcal{C}$.
\end{lem}
\begin{proof}
We will use the group law on $\mathcal{C}$ and hence we fix a
point $p_0\in \mathcal{C}$ serving as identity. Then we choose a
point $p_2\in \mathcal{C}$ such that $2p_2=p_0$ and a generic
point $p_1 \in \mathcal{C}$. Consider the line joining $p_1$ and $p_2$ and let
$p_3$ be the third intersection point with $\mathcal{C}$. Joining
$p_3$ and $p_0$ and taking the third intersection we get the point
$p_1+p_2$. Now join $p_1+p_2$ and $p_2$ and let $p_4$ be the third
intersection point. Then joining $p_4$ and $p_0$ we get
$p_1+2p_2=p_1$. Hence, the points
$p_0,p_1,p_2,p_1+p_2,p_3,p_4\in\mathcal{C}$ are a star
configuration.
\end{proof}

We now consider the general situation:
\begin{thm}\label{l=4prop}
Let $\mathcal{C}\subset\PP^2$ be a generic degree $d\geq3$ plane curve.
Then there exists a star configuration $\mathbb{X}(4)\subset
\mathcal{C}$.
\end{thm}

\begin{proof}
The case $d=3$ is treated in Lemma \ref{grouplawlem}, while for
$d>3$ we will use Lemma \ref{tgspaceideallem} and the notation
introduced therein.
For $d=4$, we produced an explicit example in Example
\ref{l4d4ex} where $I_4=S_4$, and this
is enough to conclude by semi-continuity.

To deal with the general case $d\geq 5$, we use the structure of
the coordinate ring of a star configuration. Given four generic
linear forms $L_1,L_2,L_3, L_4$ we consider the star configuration
they define, say
$\mathbb{X}=\left\{p_{i,j} : 1\leq i <j \leq 4\right\}.$
The coordinate ring of $\mathbb{X}$ is
\[A={S\over \left(\hat L_1,\dots,\hat L_4\right)}.\]

When $d\geq 2$, $\dim_\CC A_d=6$.  To prove that $I_d = S_d$,  we want to find linear forms
$N_1,\ldots,N_6$ such that the six forms $N_iQ_i$ are linearly
independent in $A_4$.
To check whether elements in $A$ are linearly independent it is
enough to consider their evaluations at the points $p_{i,j}$.
Consider the following evaluation matrix
\begin{equation}\label{evalmatrix}
\begin{array}{c|cccc}
        & Q_1       & Q_2       & Q_3      & Q_4 \\
\hline \\

p_{1,2} & M_2L_3L_4 & M_1L_3L_4 & 0        & 0 \\

p_{1,3} & M_3L_2L_4 & 0         & M_1L_2L_4& 0  \\

p_{1,4} & M_4L_2L_3 & 0         & 0        & M_1L_2L_3\\

p_{2,3} & 0         & M_3L_1L_4 & M_2L_1L_4& 0 \\

p_{2,4} & 0         & M_4L_1L_3 & 0        & M_2L_1L_3\\

p_{3,4} & 0         & 0         & M_4L_1L_2  & M_3L_1L_2
\end{array}
\end{equation}
obtained by evaluating each $Q_i$ at the points of the
configuration where we denote, by abuse of notation,
$M_hL_iL_j(p_{m,n})$ with $M_hL_iL_j$.   For example,
since $Q_2 = M_1L_3L_4+M_3L_1L_4 + M_4L_1L_3$, $Q_2(p_{2,3})
= M_3L_1L_4(p_{2,3})$ because $L_3$ vanishes at $p_{2,3}$,
while $Q_2(p_{1,4}) = 0$ since $L_4$ and $L_1$ vanish at $p_{1,4}$.
Now choose the forms
$M_i$ with $\deg M_i = d-3 \geq 2$ so that
\[M_2(p_{2,3})=M_3(p_{3,4})= M_4(P_{2,4})=0 \]
and no other vanishing occurs at the points of the star
configuration. Notice that this is possible as $\dim_\CC A_d=6$ for
$d\geq 2$. The matrix \eqref{evalmatrix} can be represented as
\begin{equation}\label{evalmatrixstar}
\begin{array}{c|cccc}
        & Q_1       & Q_2       & Q_3      & Q_4 \\
\hline \\

p_{1,2} & * & * & 0 & 0 \\

p_{1,3} & * & 0 & * & 0  \\

p_{1,4} & * & 0 & 0 & *  \\

p_{2,3} & 0 & * & 0 & 0  \\

p_{2,4} & 0 & 0 & 0 & *  \\

p_{3,4} & 0 & 0 & * & 0
\end{array}
\end{equation}
where $*$ denotes a non-zero scalar. As this matrix has rank four,
then the forms $Q_i$ are linearly independent in $A$.

To represent $N_iQ_i$ in $A$ it is enough to multiply the $j$-th
element of the $i$-th column of \eqref{evalmatrix} by the
evaluation of $N_i$ at $p_{j,i}$. Hence, if we choose the forms
$N_i$ such that they do not vanish at any point of $\mathbb{X}$,
the evaluation matrix of the $N_iQ_i$ has the same non-zero pattern
as \eqref{evalmatrixstar}. Thus the forms $N_1Q_1,\ldots,N_4Q_4$
are linearly independent in $A$.

To complete the proof we need two more forms, and we choose
$L_1Q_3$ and $L_1Q_2$ whose evaluation matrix at the points
$p_{i,j}$ is
\begin{equation}\label{evalmatrix2more}
\begin{array}{c|cccccc}
    & p_{1,2}     & p_{1,3}     & p_{1,4}     & p_{2,3}     & p_{2,4} & p_{3,4} \\ \hline \\
L_1Q_3 & 0 & 0 & 0 & 0           & 0 & * \\

L_1Q_2 & 0 & 0           & 0           & * & 0 & 0
\end{array}
\end{equation}
These rows are linearly independent with the columns of
\eqref{evalmatrixstar}, completing the proof.
\end{proof}

\section{The $l=5$ case}

It remains to consider the case that $l=5$. The case $l=5$ and
$d=4$ was classically studied.  A quartic containing a star
configuration of ten points, that is, a $\XX(5)$.   For more
details, see \cite{luroth,morley}, and  for a modern treatment we
refer to \cite{ottaviani sernesi}. In particular, we require the
following property of L\"uroth quartics:

\begin{thmNONUMBER}[Theorem 11.4 of \cite{ottaviani sernesi}]
{\it L\"{u}roth quartics form a hypersurface of degree $54$ in the
space of plane quartics.}
\end{thmNONUMBER}

\begin{cor}  \label{luroth}
Question \ref{mainproblem} has a negative
answer for $(d,l) = (4,5)$.
\end{cor}

\begin{proof}
We notice that the previous theorem is enough to give a negative
answer to our question for $l=5$ and $d=4$. In fact, a
generic plane quartic is not a L\"{u}roth quartic, and hence,
no star configuration $\mathbb{X}(5)$ can be found on it.
\end{proof}

For the remaining values $d >4$, we give an affirmative
answer to Question \ref{mainproblem}:

\begin{thm}\label{l=5}
The generic degree $d>4$ plane curve contains a star configuration
$\mathbb{X}(5)$.
\end{thm}

\begin{proof}
For the case $d=5$ we produce an explicit example and then we
conclude by semi-continuity.   For the general case, we
produce a proof similar to Theorem \ref{l=4prop}.

For the $d=5$ case, let our linear forms be given by
\[L_1 = x, ~~ L_2 = y, L_3 = z, ~~ L_4 = x+y+z, ~~\text{and}~~ L_5 = 2x-3y+5z.\]
To construct the polynomials $Q_1,\ldots,Q_5$, we make use of the following forms:
\[M_1 = x+y-z,~~ M_2 = -x+2y+2z,~~ M_3 = 2x-y-z,~~ M_4 = 3x+y-z, ~~\text{and}~~M_5 = 4x-4y+3z.\]
Form the ideal $I$ as in Lemma \ref{tgspaceideallem}.  Using
CoCoA \cite{C} to compute $I_5$, we find that $I_5 = S_5$.

For the $d\geq 6$ case, we use the notation of Lemma \ref{tgspaceideallem}, and
fix linear forms $L_1,\ldots,L_5$, and use $p_{i,j}$ to denote the ten
points of $\XX(5)$ defined by the linear forms $L_i$.   We construct the evaluation table

\begin{equation}\label{evalmatrixl=5}
\begin{array}{c|ccccc}
    & Q_1 & Q _2 & Q_3 & Q_4 & Q_5 \\ \hline \\

p_{1,2} & M_2 L_3 L_4 L_5 & M_1 L_3 L_4 L_5 & 0 & 0 & 0 \\

p_{1,3} & M_3L_2L_4L_5 & 0 & M_1 L_2L_4L_5 & 0 & 0 \\

p_{1,4} & M_4 L_2L_3L_5 &  0 & 0 & M_1L_2L_3L_5 & 0 \\

p_{1,5} & M_5 L_2 L_3 L_4 & 0 & 0 & 0 & M_1 L_2 L_3 L_4 \\

p_{2,3} & 0 & M_3 L_1L_4L_5 & M_2 L_1L_4L_5 & 0 & 0 \\

p_{2,4} & 0 & M_4 L_1L_3L_5 & 0 & M_2 L_1L_3L_5 & 0 \\

p_{2,5} & 0 & M_5 L_1L_3L_4 & 0 & 0 & M_2 L_1L_3L_4 \\

p_{3,4} & 0 & 0 & M_4 L_1L_2L_5 & M_3 L_1 L_2 L_5 & 0 \\

p_{3,5} & 0 & 0 & M_5L_1L_2L_4 & 0 & M_3 L_1L_2L_4 \\

p_{4,5} & 0 & 0 & 0 & M_5 L_1L_2L_3 & M_4 L_1L_2L_3
\end{array}
\end{equation}
obtained by evaluating each $Q_i$ at the points of the
configuration where we denote, by abuse of notation,
$M_hL_iL_jL_r(p_{m,n})$ with $M_hL_iL_jL_r$.

We work in the coordinate ring of the star configuration
$A={S\over \left(\hat L_1,\dots,\hat L_5\right)}$. Now consider
\begin{equation}\label{formsl=5}
 L_5Q_1,L_2Q_1, L_1Q_2,L_3Q_2, L_2Q_3,L_4Q_3,
L_3Q_4,L_5Q_4, L_4Q_5,L_1Q_5
\end{equation}
and we want to show that they are linearly independent in $A$ for
a generic choice of the forms $M_i$, $\deg M_i=d-4$. Again, it is
enough to show this for a special choice of forms. We choose the
forms $M_i$ in such a way that
\[M_1(p_{1,5})=M_4(p_{3,4})=M_5(p_{4,5})=0,\]
\[M_2(p_{1,2})=M_2(p_{2,5})=M_3(p_{1,3})=M_3(p_{2,3})=0,\]
and the following are non-zero
\[M_1(p_{1,2}),M_1(p_{1,3}),M_1(p_{1,4}),~~~~~M_2(p_{2,3}),M_2(p_{2,4}),~~~~~~M_3(p_{3,4}),M_1(p_{3,5}),\]
\[M_4(p_{1,4}),M_4(p_{2,4}),M_4(p_{4,5}),~~~~~M_5(p_{1,5}),M_5(p_{2,5}),M_5(p_{3,5}).\]
Notice that these conditions can be satisfied  because $\deg
M_i\geq 2$ and $\dim_\CC A_2=6, \dim_\CC A_e=10$ for $e\geq 3$.
Then, evaluating the forms in \eqref{formsl=5} at the points
$p_{i,j}$, we obtain the matrix

\begin{equation}\label{starmatrixl=5}
\begin{array}{c|cccccccccc}
       & L_5Q_1& L_2Q_1& L_1Q _2& L_3Q_2& L_2Q_3& L_4Q_3& L_3Q_4& L_5Q_4& L_4Q_5& L_1Q_5\\ \hline \\

p_{1,2}& 0     & 0     & 0      & *      & 0    & 0     & 0     & 0     & 0     & 0      \\

p_{1,3}& 0     & 0     & 0      & 0      & *    & *     & 0     & 0     & 0     & 0      \\

p_{1,4}& *     & *     & 0      & 0      & 0    & 0     & *     & *     & 0     & 0      \\

p_{1,5}& 0     & *     & 0      & 0      & 0    & 0     & 0     & 0     & 0     & 0      \\

p_{2,3}& 0     & 0     & 0      & 0      & 0    & *     & 0     & 0     & 0     & 0      \\

p_{2,4}& 0     & 0     & *      & *      & 0    & 0     & *     & *     & 0     & 0       \\

p_{2,5}& 0     & 0     & *      & *      & 0    & 0     & 0     & 0     & 0     & 0      \\

p_{3,4}& 0     & 0     & 0      & 0      & 0    & 0     & 0     & *     & 0     & 0      \\

p_{3,5}& 0     & 0     & 0      & 0      & *    & *     & 0     & 0     & *     & *      \\

p_{4,5}& 0     & 0     & 0      & 0      & 0    & 0     & 0     & 0     & 0     & *      \\
\end{array}
\end{equation}
where $*$ denotes a non-zero scalar. One can
verify that \eqref{starmatrixl=5} has rank ten and hence the
forms in \eqref{formsl=5} are linearly independent in $A$ and the
thesis follows.
\end{proof}

We can now prove our main theorem.

\begin{thm}\label{maintheorem}
Let $l \geq 2$.  Then the generic degree $d$ plane curve
contains a star configuration $\XX(l)$ if and only if
\begin{enumerate}
\item[$(i)$] $l=2$ and $d \geq 1$, or
\item[$(ii)$] $l=3$ and $d \geq 2$, or
\item[$(iii)$] $l=4$ and $d \geq 3$, or
\item[$(iv)$] $l = 5$ and $d \geq 5$.
\end{enumerate}
\end{thm}

\begin{proof}
Combine Lemmas \ref{trivialcase} and \ref{grouplawlem},
Corollaries  \ref{d<l-1rem} and \ref{luroth}, and
Theorems \ref{asymptoticresult}, \ref{l=4prop}, and \ref{l=5}.
\end{proof}

%%%%%%%%%%%%%%%%%%%%%%%%%%%%%%%%%%%%%%%%%%%%%%%%%%%%%%%%%%%%%%%%%%%%%%%%%%%%%%

\end{document}